\documentclass{elsarticle}
\usepackage{txfonts}
\usepackage{graphicx}
\usepackage{tikz}

\usetikzlibrary{positioning}

\newtheorem{thm}{Theorem}

\begin{document}

\title{On the F-index and F-coindex of the line graphs of the subdivision graphs}

\author[man]{Ruhul Amin}
\ead{aminruhul80@gmail.com}
\author[man]{Sk. Md. Abu Nayeem\corref{cor1}}
\ead{nayeemsma@gmail.com}

\address[man]{Department of Mathematics, Aliah University, New Town, Kolkata -- 700 156, India.}
\cortext[cor1]{Corresponding Author.}

\begin{abstract}
The aim of this work is to investigate the F-index and F-coindex of the line graphs of the cycle graphs, star graphs, tadpole graphs, wheel graphs and ladder graphs using the subdivision concepts. F-index of the line graph of subdivision graph of square grid graph, 2D-lattice, nanotube and nanotorus of $TUC_{4}C_{8}[p, q]$ are also investigated here.

\end{abstract}


\maketitle

\section{Introduction}
Topological indices are numerical quantities associated with different graph parameters. These are used to correlate chemical structure of molecular graphs with various physical properties, chemical reactivities and biological activities. By molecular graph we mean a simple graph, representing the carbon-atom skeleton of an organic molecule. Let $G$ be a simple graph with vertex set $V(G)$, edge set $E(G)$ and $d(u)$ denotes the degree of a vertex $u$ in $G$. An edge between vertices $u$ and $v$ is denoted by $uv$. The set of vertices which are adjacent to a vertex $u$ is denoted by $N(u)$ and is known as neighbourhood of $u$. Clearly, $d(u)=|N(u)|$. Degree based topological indices have been subject to study since the introduction of Randi\'{c} index in 1975. Although the first degree-based topological indices are the Zagreb indices \cite{gut72}, these were initially intended for the study of total $\varphi$-electron energy~\cite{gut13cca} and were included among the topological indices much later. The first and second Zagreb indices are respectively defined as $$M_{1}(G) = \displaystyle\sum_{u \in V(G)}d^2(u) = \displaystyle\sum_{uv \in E(G)}[d(u) + d(v)],$$  and $$M_{2}(G) = \displaystyle\sum_{uv \in E(G)}d(u)d(v).$$

In the same paper \cite{gut72} where Zagreb indices were introduced, Gutman and Trinajsti\'{c} indicated that another term of the form $\displaystyle\sum_{u \in V(G)}d^3(u)$ influences the total $\varphi$-electron energy. But this remained unstudied by the researchers for a long time, except for a few occasions \cite{hu05, li04, li05} until the publication of an article by Furtula and Gutman in 2015 and so they named it ``forgotten topological index'' or F-index in short~\cite{fur15}. Thus, F-index of a graph $G$ is defined as $$F(G) = \displaystyle\sum_{u \in V(G)}d^3(u) = \displaystyle\sum_{uv \in E(G)}[d^2(u) + d^2(v)].$$

F-index for different graph operations has been studied De et al.~\cite{de16a}. Extremal trees with respect to F-index have been studied by Abdo et al.~\cite{abo15}.

While considering contribution of non adjacent pair of vertices in computing the weighted Wiener polynomials of certain composite graphs, Do\v{s}lic~\cite{dos08} introduced quantities named as Zagreb co-indices. Thus the first and second Zagreb co-indices are respectively defined as $$\overline M_{1}(G) = \displaystyle\sum_{uv \notin E(G)}[d(u) + d(v)]$$ and $$\overline M_{2}(G) = \displaystyle\sum_{uv \notin E(G)}[d(u)d(v)].$$

In a similar manner, F-coindex is defined as $$\overline{F}(G)=\sum_{uv\not\in E(G)}[d^2(u)+d^2(v)].$$ De et al. \cite{de16b} have shown that F-coindex can predict the octanol water partition coefficients of molecular structures very efficiently. They have also studied the F-coindex of graph operations. Trees with minimum F-coindex have been found by Amin and Nayeem \cite{ami16}.

Complement of the graph $G$ is a simple graph $\overline{G}$ with same vertex set $V$ and there is an edge between the vertex $u, v$ in $\overline{G}$ if and only if there is no edge between $u, v$ in $G$. It is evident that the F-coindex of $G$ is not same as F-index of $\overline{G}$, because during computation of F-coindex, degrees of the vertices are considered as in $G$.

\medskip

The subdivision graph $S(G)$ is the graph obtained from $G$ by replacing each of its edge by a path of length $2$, or equivalently, by inserting an additional vertex into each edge of $G$. The line graph of the graph $G$, denoted by $L(G)$, is the simple graph whose vertices are the edges of $G$, with $ef$ belong to $E(L(G))$ when $e$ and $f$ are incident to a common vertex in $G$. We use $P_{n}, C_{n}, S_{n}$ and $W_{n}$ to denote path, cycle, star and wheel graphs on $n$ vertices respectively. The tadpole graph $T_{n,k}$ is the graph obtained by joining a cycle graph $C_{n}$ to a path of length $k$. The ladder graph $L_{n}$ is given by $L_{n} = P_{2}\Box P_{n}$, the cartesian product of $P_{2}$ with $P_{n}$. The graph obtained via this definition has the advantage of looking like a ladder, having two rails and $n$ rungs between them. In a similar manner, the square grid graph $G(m, n)$ is defined by $G(m, n) = P(m)\Box P(n)$. Thus $L(n) = G(2, n)$.

\begin{figure}
\begin{center}
\includegraphics[width=\textwidth]{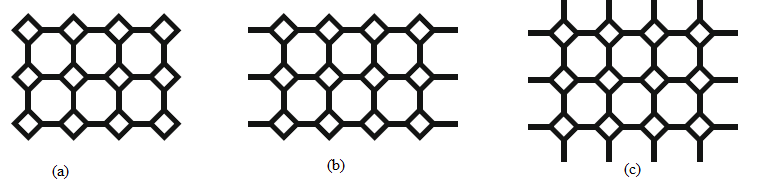}
\caption{\label{nano}(a) 2D-lattice of $TUC_4C_8[4$,3]. (b)$TUC_4C_8[4,3]$ nanotube. (c)$TUC_{4}C_{8}[4, 3]$ nanotorus. }
\end{center}
\end{figure}

Study on different topological indices for line graphs of subdivision graphs have been done by many researchers in the recent past. Ranjini et al. \cite{ran11} computed Zagreb index and Zagreb co-index of the line graphs of the subdivision graphs of $T(n, k), L(n)$ and $W(n + 1)$. Su and Xu \cite{su15} generalized the idea of Ranjini et al. and presented the Schur-bound for the general sum-connectivity co-index. Nadeem et al. \cite{nad15} studied the $ABC_4$ and $GA_5$ indices of the line graph of tadpole, wheel and ladder graphs using the notion of subdivision. Recently Nadeem et. al \cite{nad16} obtained expressions for certain topological indices for the line graph of subdivision graphs of 2D-lattice, nanotube and nanotorus of $TUC_{4}C_{8}[p, q]$, where $p$ and $q$ denote the number of squares in a row and number of rows of squares, respectively in 2D-lattice, nanotube and nanotorus. 2D-lattice, nanotube and nanotorus of $TUC_{4}C_{8}[p, q]$ are depicted in Figure \ref{nano}. The order and size of 2D-lattice of $TUC_{4}C_{8}[p, q]$ are $4pq$ and $6pq - p - q$ respectively. The order and size of nanotube of $TUC_{4}C_{8}[p, q]$ are $4pq$ and $6pq - p$ respectively. Again the order and size of nanotorus of $TUC_{4}C_{8}[p, q]$ are $4pq$ and $6pq$ respectively. Hosamani and Zafar~\cite{hos} have studied some more topological indices of the line graphs of the subdivision graphs of the above mentioned nano-structures. Here we have computed the F-index and F-coindex of the line graphs of subdivision graph of $C(n), S(n), T(n, k), L(n)$ and $W(n + 1)$ in Section 2. We also study the line graph of subdivision graphs of 2D-lattice, nanotube and nanotorus of $TUC_{4}C_{8}[p, q]$ and calculate the F-index of the line graph of subdivision of 2D-lattice, nanotube and nanotorus of $TUC_{4}C_{8}[p, q]$ in Section 3.

\section{F-index and F-coindex of some graph using subdivision concept }

\begin{thm}\label{17}
Let $G$ be the line graph of the subdivision graph of the cycle $C_{n}$ with $n$ vertices. Then F-index of $G$ is $ F(G) = 16n$.
\end{thm}

\begin{figure}
\begin{center}
\includegraphics[width=0.6\textwidth]{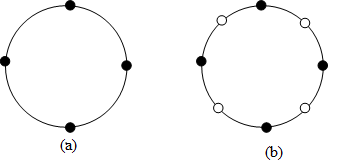}
\caption{\label{cycle}(a) Cycle graph $C_4$. (b) Subdivision graph and line graph of subdivision graph of $C_4$.}
\end{center}
\end{figure}

\noindent\textit{Proof.}
The line graph of the subdivision graph of the cycle $C_{n}$ is $C_{2n}$, i.e., $G=C_{2n}$ (See Figure \ref{cycle}). Hence
\begin{eqnarray*}
F(G) &=& \sum_{i=1}^{2n} 2^3\\
  &=&2n.2^3\\
&=& 16n.
\end{eqnarray*}
\qed

\begin{thm}\label{18}
Let $G$ be the line graph of the subdivision graph of the cycle $C_{n}$ with $n$ vertices. Then F-coindex of $G$ is $\overline{F}(G) = 16n^2 - 16n$.
\end{thm}
\noindent\textit{Proof.}
As before, $G = C_{2n}$. Each vertex of $G$ is of degree two. Also each vertex of $G$ is non-adjacent to $2n - 2$ two degree vertex. Hence
$$\sum_{u\in V(G)}\sum_{v\not \in N(u)} [d^2(u) + d^2(v)] = 2n(2n - 2)(2^2 + 2^2).$$

Since one edge is shared by a pair of vertices, $$\overline{F}(G) = \frac{1}{2}2n(2n - 2)(2^2 + 2^2) = 16n^2 - 16n.$$
\qed

\begin{figure}
\begin{center}
\includegraphics[width=0.5\textwidth]{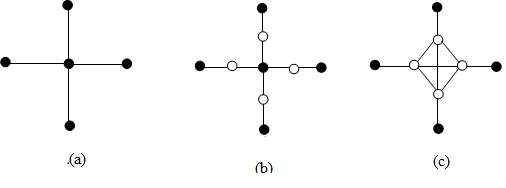}
\caption{\label{star}(a) Star graph $S_4$. (b) Subdivision of $S_4$. (c) Line graph of subdivision graph of $S_4$.}
\end{center}
\end{figure}

\begin{thm}\label{19}
Let $G$ be the line graph of the subdivision graph of the star $S_n$ with $n$ vertices (See Figure \ref{star}). Then F-index of $G$ is $F(G) = n(n- 1)(n^2 - 3n + 3)$.
\end{thm}
\noindent\textit{Proof.}
The line graph of the subdivision graph of the star contains $n -1$ vertices of degree $n - 1$ and $n - 1$ vertices of degree one. Hence
\begin{eqnarray*}
F(G) &=&\sum_{v\in V(G)} d^3(v)\\
&=&(n -1)(n -1)^3 + (n - 1)1^3\\
&=&(n - 1)(n^3 - 3n^2 + 3n - 1 + 1)\\
&=&n(n - 1)(n^2 - 3n + 3).
\end{eqnarray*}
\qed

\begin{thm}\label{20}
Let $G$ be the line graph of the subdivision graph of the star $S_n$ with $n$ vertices. Then F-coindex of $G$ is $\overline{F}(G) = (n -1)(n- 2)(n^2 - 2n + 3)$.
\end{thm}
\noindent\textit{Proof.}
Each one degree vertex, i.e., pendant vertex in $L(S(S_n))$ is non-adjacent with $n - 2$ one degree vertex and $n - 2$ vertices of degree $n - 1$. All non-pendant vertices are adjacent to each other (See Figure \ref{star}). Hence
\begin{eqnarray*}
\overline{F}(G) &=&\frac{1}{2}(1^2 + 1^2)(n - 2)(n - 1) + (n - 1)(n - 2)\{1^2 + (n - 1)^2\}\\
&=&(n - 1)(n - 2)(n^2 - 2n + 3).
\end{eqnarray*}
\qed

\begin{figure}
\begin{center}
\includegraphics[width=0.75\textwidth]{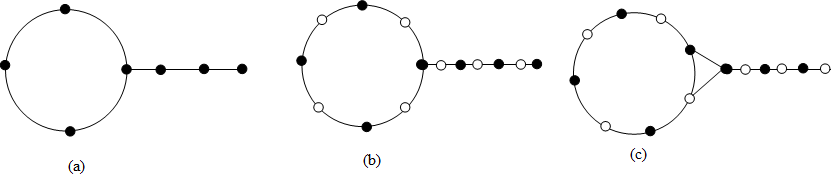}
\caption{\label{tadpole}(a) Tadpole graph $T_{4,3}$. (b) Subdivision of $T_{4,3}$. (c) Line graph of subdivision graph of $T_{4,3}$.}
\end{center}
\end{figure}

\begin{thm}\label{12}
For the line graph of the subdivision graph of a tadpole graph, $F(L(S(T_{n,k}))) = 16n + 16k + 50$.
\end{thm}
\noindent\textit{Proof.}
The subdivision graph $S(T_{n,k})$ contains $2(n + k)$ edges, so that the line graph contains $ 2(n + k)$ vertices, out of which three vertices are of degree $3$, one vertex is of degree $1$ and the remaining $2n + 2k - 4$ vertices are of degree $2$ (See Figure \ref{tadpole}).

Hence,
\begin{eqnarray*}
  F(L(S(T_{n,k})))&=&3.3^3+1.1^3+(2n+2k-4).2^3\\
&=&16n + 16k + 50.
\end{eqnarray*}
\qed

\begin{thm}\label{13}
The F-coindex of the line graph $L(S(T_{n, k}))$ is $16(n + k)^2 - 57$.
\end{thm}
\noindent\textit{Proof.}
The line graph $ L(S(T_{n, k}))$ contains a subgraph $P_{2k - 1}$. $L(S(T_{n, k}))$ contains only $3 $ vertices of degree $3$. Let $u,v_1,v_2$ be those vertices of degree $3$, among which the vertex $u$ is attached to the path $P_{2k - 1}$. The vertex $u$ is not adjacent to $2k - 3$ vertices of degree two and the pendant vertex among the $2k-1$ vertices on $P_{2k-1}$. The neighbor of $u$ on $P_{2k-1}$ is not adjacent with $2k - 4$ vertices of degree $2$ and the pendant vertex among the $2k-1$ vertices on $P_{2k-1}$, and so on. The vertices $v_1$ and $v_2$ are not adjacent with $2n + 2k - 5$ vertices of degree $2$ and also with the pendent vertex. The vertex $u$ is not adjacent with $2n - 2$ vertices of degree $2$ in $L[S(C_{n}) + e]$, where $e$ is the edge adjacent to $S(C_{n})$. Also, $2n - 2$ vertices of degree $2$ in $L[S(C_{n}) + e]$ are not adjacent with any of the vertices on the path $P_{2k - 1}$. Out of these $2n -2$ vertices in $L[S(C_{n}) + e]$, $2n - 4$ vertices of degree $2$  are non adjacent with $2n - 5$ vertices of degree $2$ and the remaining $2$ vertices which are adjacent to $v_1$ and $v2$ have $2n - 4$ non adjacent vertices in $L[S(C_{n}) + e]$ of degree $2$. Hence,
\begin{eqnarray*}
\overline{F}(L(S(T_{n ,k})))&=&(2k - 2)(1^2 + 2^2) + (2n - 2)(1^2 + 2^2) + 3 (2k + 2n - 5)(2^2 + 3^2)\\
&& + (2n - 2)(2k - 2)(2^2 + 2^2) + 3 (1^2 + 3^2)\\
&& + \{(2^2 + 2^2) + 2 (2^2 + 2^2) + \ldots + (2k - 4)(2^2 + 2^2)\}\\
&& +  \{(2^2 + 2^2) + 2 (2^2 + 2^2) + \ldots + (2n - 4)(2^2 + 2^2)\}\\
&=& 5(2k - 2) + 5 (2n - 2) + 39 (2n + 2k - 5) + 8 (2n - 2)(2k - 2)\\
 &&+ 30 + 4 (2k - 4 )(2k - 3) + 4 (2n - 4)(2n - 3)\\
&=& 16n^2 + 16k^2 + 32 n k - 57\\
&=& 16 (n + k)^2 - 57.
\end{eqnarray*}
\qed

\begin{thm}\label{13}
The F-index of $L[S(W_{n + 1})]$ is $n (n^3 + 81)$.
\end{thm}
\begin{figure}
\begin{center}
\includegraphics[width=0.75\textwidth]{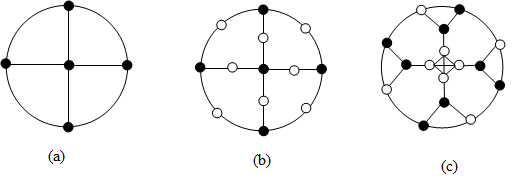}
\caption{\label{wheel}(a) Wheel graph. (b) Subdivision of wheel graph. (c) Line graph of subdivision graph of wheel graph}
\end{center}
\end{figure}
\noindent\textit{Proof.}
The line graph $L[S(W_{n + 1})]$ contains $4n$ vertices. Out of these $4n$ vertices, $3n$ vertices are of degree $3$ and $n$ vertices are of degree $n$ (See Figure \ref{wheel}).

Hence

\[F(G)=3n.3^3+n.n^3=n(n^3+81).\]
\qed

\begin{thm}\label{14}
The F-coindex of $L[S(W_{n + 1})]$ is $2n (n^3 + 56n - 56)$.
\end{thm}
\noindent\textit{Proof.}
As seen in Theorem \ref{13}, $L(S(W_{n+1}))$ has total $4n$ vertices, out of which $3n$ vertices are of degree $3$ and $n$ vertices are of degree $n$. Out of the $3n$ three degree vertices, $2n$ vertices are not adjacent to the $3n - 4$ vertices of degree $3$ and $n$ vertices of degree $n$. Hence their contribution in the F-coindex is
\begin{eqnarray*}
&&2n(3n - 4)(3^2 + 3^2) + 2n.n(3^2 + n^2)\\
&=&108n^2 - 144n + 18n^2 + 2n^3.
\end{eqnarray*}

The remaining $n$ vertices of degree $3$ are non adjacent with $n - 1$ vertices of degree $n$ and $3n - 3$ vertices of degree $3$. Hence their contribution to the F-coindex is
\begin{eqnarray*}
&&n(n - 1)(3^2 + n^2) + n(3n - 3)(3^2 + 3^2)\\
&=& (n^2 - n)(9 + n^2) + 18n(3n - 3)\\
&=& n^4 - n^3 + 63n^2 - 63n.
\end{eqnarray*}

Also each of the $n$ vertices of degree $n$ are non adjacent with $3n - 1$ vertices of degree $3$. So, their contribution to F-coindex is
\begin{eqnarray*}
&&n(3n - 1)(3^2 + n^2)\\
&=& (3n^2 - n)(9 + n^2)\\
&=& 3n^4 + 27n^2 - 9n - n^3.
\end{eqnarray*}

Hence,
 \begin{eqnarray*}
\sum_{u\in V(G)}\sum_{u\not\in N(u)}[(d(u))^2 + (d(v))^2] &= &2n^3 + 18n + 108n^2 - 144n + n^4 - n^3 + 63n^2 - 63n\\
&&+3n^4 + 27n^2 - n^3 - 9n\\
&=& 4n(n^3 + 56n - 56).
\end{eqnarray*}

Since each edge is shared by a pair of vertices, $$\overline{F}(L(S(W_{n + 1}))) = 2n(n^3 + 56n - 56).$$
\qed

\begin{thm}\label{15}
The F-index of $L(S(L_{n}))$ is $162n - 260$.
\end{thm}
\begin{figure}
\begin{center}
\includegraphics[width=\textwidth]{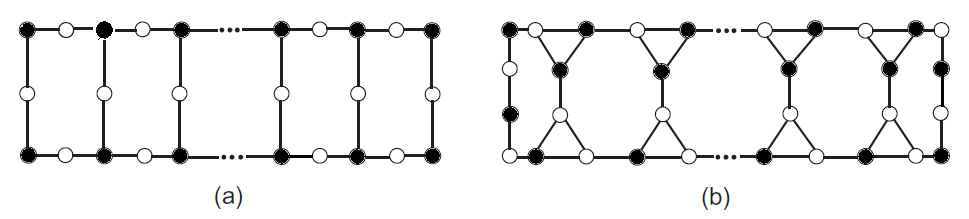}
\caption{\label{ladder}(a) Subdivision of lader graph. (b)Line graph of subdivision graph of lader graph.}
\end{center}
\end{figure}
\noindent\textit{Proof.}
$L(S(L_{n}))$ contains total $6n - 4$ vertices (See Figure \ref{ladder}). Out of these $6n - 4$ vertices, $8$ vertices are of degree $2$ and the remaining $6n - 12$ vertices are of degree $3$. Hence,
\begin{eqnarray*}
F(L(S(L_{n})))=8.2^3+(6n-12).3^3=162n-260.
\end{eqnarray*}
\qed

\begin{thm}\label{16}
The F-coindex of the line graph $L(S(L_{n}))$ is $324n^2 - 832n + 532$.
\end{thm}
\noindent\textit{Proof.}
As stated earlier, $L(S(L_n))$ has $6n-4$ vertices in total, among which 8 vertices are of degree two and all other vertices are of degree three. Out of the $6n-12$ three degree vertices, the 4 vertices which are adjacent to the corner vertices are not adjacent to $6n-14$ three degree vertices and 7 two degree vertices. hence their contribution in the F-coindex is $$4(6n-14)(3^2+3^2)+20(2^2+2^2).$$

Contribution of the rest of the three degree vertices, each of which is not adjacent to $6n-15$ three degree vertices and 8 two degree vertices, is $$(6n-16)(6n-15)(3^2+3^2)+8(6n-16)(2^2+3^2).$$

Each of the four corner vertices of degree two is not adjacent to $6n-13$ three degree vertices and 6 two degree vertices. Hence their contribution in the F-coindex is $$4(6n - 13) (2^2 + 3^2) + 24 (2^2 + 2^2).$$

Contribution of the rest of the two degree vertices is $$4(6n - 12) (2^2 + 3^2) + 20 (2^2 + 2^2).$$

Hence,
 \begin{eqnarray*}
\sum_{u\in V(G)}\sum_{v\not\in N(u)}[(d(u))^2 + (d(v))^2] &= &18 (6n - 16) (6n - 15) + 104 (6n - 16) + 72 (6n - 14)\\
&&+ 364 + 52 (6n - 12)+ 160 + 52 (6n - 13) + 192\\
&=& 648n^2 - 1668n  + 1064.
\end{eqnarray*}

Since one edge is shared by a pair of vertices, $$\overline{F}(L(S(L_{n}))) = \frac{1}{2}(648n^2 - 1668n + 1064) = 324n^2 - 832n + 532.$$
\qed

\section{F-index of some nanostructure}
\begin{figure}
\begin{center}
\includegraphics[width=\textwidth]{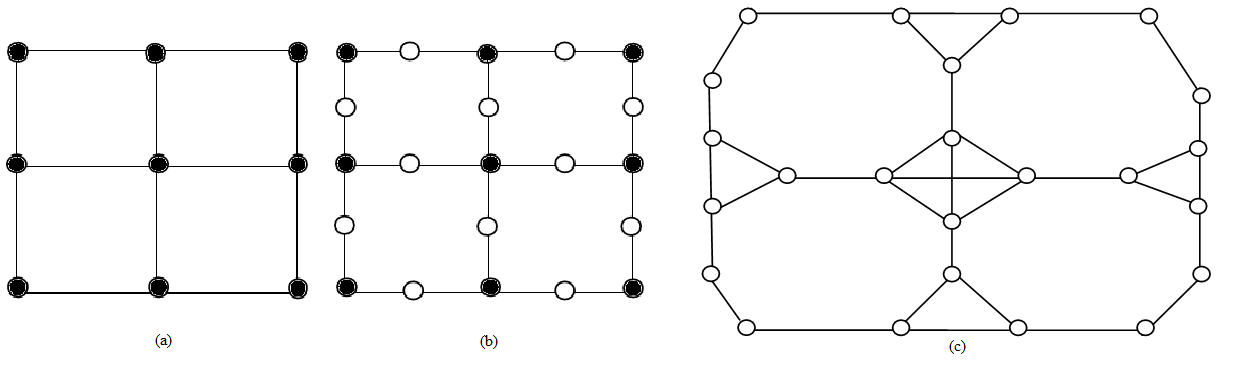}
\caption{\label{grid}(a)Square grid graph. (b) Line graph of square grid graph. (c) Subdivision graph of line graph of square grid graph.}
\end{center}
\end{figure}

\begin{thm}\label{21}
Let $G$ be a line graph of the subdivision graph of a square grid graph with $mn$ vertices (See Figure \ref{grid}). Then F-index of $G$ is $F(G) = 256mn - 350m - 350n + 440$.
\end{thm}
\noindent\textit{Proof.}
The line graph of the subdivision graph of a square grid graph with $mn$ vertices, i.e., $(m - 1)(n - 1)$ number of squares contains $8$ two degree vertices, $6(n - 2) + 6(m - 2)$ three degree vertices and $4(m - 2)(n - 2)$ four degree vertices. Hence,
\begin{eqnarray*}
F(G) &= &\sum_{v\in V(G)} d^3(v)\\
&=& 8.2^3 + [6(n - 2) + 6(m - 2)]3^3 + 4(m - 2)(n - 2)4^3\\
&=& 64 + 162(m + n - 4) + 256(m - 2)(n - 2)\\
&=& 256mn - 350m - 350n + 440.
\end{eqnarray*}
\qed

\begin{figure}
\begin{center}
\includegraphics[width=\textwidth]{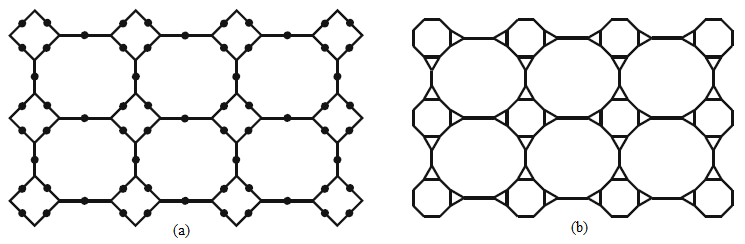}
\caption{\label{lattice}(a) Subdivision graph 2D--lattice $TUC_4C_8[4,3]$. (b) Subdivision graph of line graph of 2D--lattice.}
\end{center}
\end{figure}

\begin{thm}\label{22}
Let $G$ be a line graph of the subdivision graph of a 2D-lattice of $TUC_{4}C_{8}[p, q]$ (See Figure \ref{lattice}). Then F-index of $G$ is $F(G) = 324pq - 130p - 130q$.
\end{thm}
\noindent\textit{Proof.}
The line graph of the subdivision graph of a 2D-lattice of $TUC_{4}C_{8}[p, q]$ contains $4p + 4q$ two degree vertices and $12pq - 6p - 6q$ three degree vertices. Hence,

\begin{eqnarray*}
F(G) &=& \sum_{v\in V(G)} d^3(v)\\
&=& (4p + 4q)2^3 + (12pq - 6p - 6q)3^3\\
&=& 324pq - 162p - 162q + 32p + 32q\\
&=& 324pq -  130p - 130q.
\end{eqnarray*}
\qed

\begin{figure}
\begin{center}
\includegraphics[width=\textwidth]{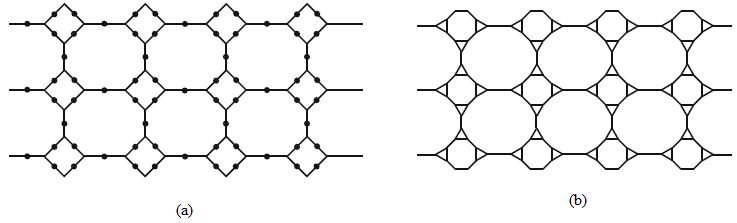}
\caption{\label{nanotube}(a) Subdivision of $TUC_{4}C_{8}$ nanotube. (b) Line graph of subdivision graph of $TUC_{4}C_{8}$ nanotube.}
\end{center}
\end{figure}

\begin{thm}\label{23}
Let $G$ be a line graph of the subdivision graph of $TUC_{4}C_{8}[p, q]$ nanotube (See Figure \ref{nanotube}). Then F-index of $G$ is $F(G) = 324pq - 130p + 2q$.
\end{thm}
\noindent\textit{Proof.}
The line graph of the subdivision graph of $TUC_{4}C_{8}[p, q]$ nanotube contains $ 2q$ pendant vertices, $4p$ two degree vertices and $12pq - 6p$ three degree vertices. Hence,
\begin{eqnarray*}
F(G) &= &\sum_{v\in V(G)} d^3(v)\\
&=& 2q.1^3 + 4p.2^3 + (12pq - 6p).3^3\\
&=& 2q + 32p + 27(12pq - 6p)\\
&=& 324pq - 130p + 2q.
\end{eqnarray*}
\qed

\begin{figure}
\begin{center}
\includegraphics[width=\textwidth]{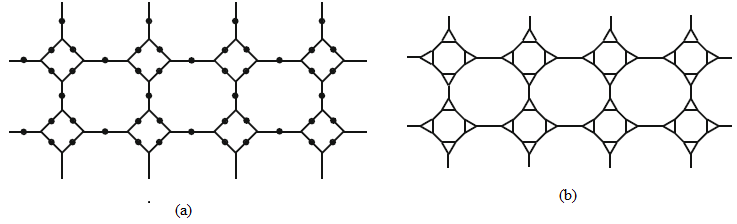}
\caption{\label{nanotorus}(a) Subdivision of $TUC_{4}C_{8}$ nanotorus. (b) Line graph of subdivision graph of $TUC_{4}C_{8}$ nanotorus.}
\end{center}
\end{figure}

\begin{thm}\label{24}
Let $G$ be a line graph of the subdivision graph of $TUC_{4}C_{8}[p, q]$ nanotorus (See Figure \ref{nanotorus}). Then F-index of $G$ is $F(G) = 324pq + 2p + 2q$.
\end{thm}
\noindent\textit{Proof.}
The line graph of the subdivision graph of $TUC_{4}C_{8}[p, q]$ nanotorus contains $2p +  2q$ pendant vertices, and $12pq$ three degree vertices. Hence,
\begin{eqnarray*}
F(G) &= &\sum_{v\in V(G)} d^3(v)\\
&=& (2p + 2q).1^3 + 12pq.3^3\\
&=& 324pq + 2p + 2q.
\end{eqnarray*}
\qed

\section*{Acknowledgement}
This work has been partially supported by University Grants Commission, India through Grant No F$1$-$17.1/2016$-$17/$MANF-$2015$-$17$-WES-$60163/$(SA-III/Website) to the first author.

\end{document}